\documentclass{amsart}
\usepackage{graphicx}
\usepackage{amssymb}
\newtheorem{thm}{Theorem}

\newtheorem{prop}[thm]{Proposition}
\theoremstyle{definition}
\newtheorem{defn}[thm]{Definition}

\theoremstyle{remark}

\numberwithin{equation}{section}
\begin{document}
\title[Continuous Functions in Variable-Exponent Lebesgue Spaces]{The Note on the Closure of Continuous Functions in Variable-Exponent Lebesgue Spaces for Multiple Variables}%
\author{Nikoloz Devdariani}%
\email[Nikoloz Devdariani]{nikoloz.devdariani275@ens.tsu.edu.ge}%

\thanks{}


%\date{}%
%\dedicatory{}%
%\commby{}%
% ----------------------------------------------------------------
\begin{abstract}
In this paper, we generalize a recently obtained result by Kopaliani and Zviadadze from the one-variable case to the several-variable case. Specifically, in terms of decreasing rearrangement, we characterize those exponents $p(\cdot) $ for which the corresponding variable-exponent Lebesgue space $L^{p(\cdot)}([0;1]^n)$ shares the property with $L^\infty([0;1]^n)$ such that the space of continuous functions $C([0;1]^n)$ forms a closed linear subspace in $L^{ p(\cdot)}([0;1]^n)$ . In particular, we derive the necessary and sufficient conditions on the decreasing rearrangement of the exponent $p(\cdot)$ for which there exists an equimeasurable exponent of $p(\cdot)$ such that the corresponding variable-exponent Lebesgue space possesses the aforementioned property.
\end{abstract}
\maketitle
% ----------------------------------------------------------------
\section{Introduction}

\par Recently, Kopaliani and Zviadadze \cite{kop-zv} characterized those exponents $p(\cdot)$ for which the corresponding variable-exponent Lebesgue space $L^{p(\cdot)}([0; 1])$ shares the property with $L^\infty([0; 1])$ that the space of continuous functions $C[0; 1]$ forms a closed linear subspace in $L^{p(\cdot)}([0; 1])$. Later, Professor Alberto Fiorenza asked the authors of the mentioned paper how a similar statement would look in the case of several variables.

\par In this paper, we provide an answer to this question. To formulate our result, let us introduce some notations. Let's define $\Omega := [0; 1]^{n}$, $I := [a_1; b_1] \times \ldots \times [a_n; b_n]$, where $0 \leq a_i < b_i \leq 1$ for $i \in {1, \ldots, n}$. We also define a measurable function $p : \Omega \to [1; \infty)$. Let $W(p)$ denote the set of all functions equimeasurable with $p(\cdot)$. Below, we will find the conditions on the function $p(\cdot)$ for which there exists $\bar{p}(\cdot) \in W(p)$ such that the space $C(\Omega)$ of continuous functions is a closed subspace in $L^{\bar p(\cdot)}(\Omega)$.

\par In particular, we prove the following theorem:

\begin{thm}
\label{thm_necess_suff_on_decreasing_rearrangement}
For the existence of $\bar{p}(\cdot) \in W(p)$ for which $C(\Omega)$ is a closed subspace in $L^{\bar{p}(\cdot)}(\Omega)$, it is necessary and sufficient that
\begin{equation}
\label{cond_necess_suffic}
\limsup_{t \to 0+}\frac{p^*(t)}{\ln(e/t)}>0.
\end{equation}
\end{thm}

\section{definitions and auxiliary results}
\par Let $\mathcal{M}$ denote the space of all equivalence classes of Lebesgue measurable real-valued functions, equipped with the topology of convergence in measure relative to each set of finite measure.

\begin{defn}
\label{defn_BFS}
A Banach subspace $X$ of $\mathcal{M}$ is referred to as a Banach function space (BFS) on $\Omega$ if the following conditions hold:

$1)$ The norm $\|f\|_{X}$ is defined for every measurable function $f$, and $f\in X$ if and only if $\|f\|_{X} < \infty$. Also, $\|f\|_{X} = 0$ if and only if $f = 0$ almost everywhere;

$2)$ $\||f|\|_{X} = \|f\|_{X}$ for all $f\in X$;

$3)$ If $0 \leq f \leq g$ almost everywhere, then $\|f\|_{X} \leq \|g\|_{X}$;

$4)$ If $0 \leq f_{n} \uparrow f$ almost everywhere, then $\|f_{n}\|_{X} \uparrow \|f\|_{X}$;

$5)$ If $E$ is a measurable subset of $\Omega$ with finite measure ($|E|<\infty$), then $\|\chi_{E}\|_{X} < \infty$, where $\chi_{E}$ is the characteristic function of $E$;

$6)$ For every measurable set $E$ with finite measure ($|E|<\infty$), there exists a constant $C_{E} < \infty$ such that $\int_{E}f(t)dt \leq C_{E}\|f\|_{X}$.
\end{defn}

\par Now, let's introduce various subspaces of a BFS $X$:
\begin{itemize}
\item A function $f$ in $X$ has an absolutely continuous norm in $X$ if $\|f\cdot\chi_{E_n}\|_{X} \to 0$ whenever ${E_n}$ is a sequence of measurable subsets of $\Omega$ such that $\chi_{E_n} \downarrow 0$ almost everywhere. The set of all such functions is denoted by $X_A$;

\item $X_B$ is the closure of the set of all bounded functions in $X$;

\item A function $f \in X$ has a continuous norm in $X$ if, for every $x \in \Omega$, $\lim_{\varepsilon \to 0+}\|f\chi_{B(x,\varepsilon)}\|_X = 0$, where $B(x,\varepsilon)$ is a ball centered at $x$ with radius $\varepsilon$. The set of all such functions is denoted by $X_C$.
\end{itemize}

\par The relationship between the concepts of $X_A$ and $X_B$ is given in \cite{BS}. Generally, the interplay among the subspaces $X_A$, $X_B$, and $X_C$ can be intricate. For instance, there exists a BFS $X$ in which ${0} = X_A \subsetneq X_C = X$ (for example see  \cite{LN}).

\par Let $\mathcal{P}$ through whale paper denotes the family of all measurable functions $p: \Omega\rightarrow[1;+\infty)$.
When $p(\cdot)\in\mathcal{P}$ we denote by $L^{p(\cdot)}(\Omega)$, the set of all measurable functions $f$ on $\Omega$ such that for some
$\lambda>0$
$$
\int_{\Omega}\left(\frac{|f(x)|}{\lambda}\right)^{p(x)}dx<\infty.
$$
This set becomes a BFS when equipped with the
norm
$$
\|f\|_{p(\cdot)}=\inf\left\{\lambda>0:\,\,
\int_{\Omega}\left(\frac{|f(x)|}{\lambda}\right)^{p(x)}dx\leq1\right\}.
$$

\par The variable exponent Lebesgue spaces $L^{p(\cdot)}(\Omega)$ and the corresponding variable exponent Sobolev spaces $W^{k,p(\cdot)}$ are of significant interest due to their applications in fluid dynamics, partial differential equations with non-standard growth conditions, calculus of variations, image processing, and more (refer to \cite{CUF,DHHR} for further details).

\par For the specific case of a particular BFS $X = L^{p(\cdot)}(\Omega)$, the relationship between this space and its subspaces, namely, $X_A$, $X_B$, and $X_C$, has been explored in \cite{ELN}. We will now present some of the key findings from that paper.

\begin{prop}[Edmunds, Lang, Nekvinda]
\label{prop_X_A_equals_X_B}
Let $p(\cdot)\in \mathcal{P}$ and set $X=L^{p(\cdot)}(\Omega)$. Then 
\par (i) $X_A=X_C$;
\par (ii) $X_B=X$ if and only if $p(\cdot)\in L^\infty(\Omega)$;
\par (iii) $X_A=X_B$ if and only if
$$
\int_0^1 c^{p^*(t)}dt<\infty, \:\:\: \text{for all}\:\:\: c>1,
$$
where $p^*$ is the decreasing rearrangement of $p(\cdot)$.
\end{prop}

\par The decreasing rearrangement of the measurable function $f$ is defined as follows:
$$
f^*(x)=\inf\{\lambda\geq0\::\:|\{|f|>\lambda\}|\leq x\}
\: \: \:  x > 0 . $$
By construction, $f^*$ is a decreasing right-continuous function. Furthermore, the functions $|f|$ and $f^*$ are
equimeasurable. 

\par If $\psi$ is an increasing convex function $\psi:[0;+\infty)\to[0;+\infty)$, such that $\psi(0)=0$,
$$
\lim_{x\to0+}(\psi(x)/x)=0,\quad \text{and} \quad \lim_{x\to+\infty}(\psi(x)/x)=+\infty,
$$
then the Orlicz space $L_\psi$ is defined as the set of all $f\in \mathcal{M}(\Omega)$ for which:
$$
||f||_{L_\psi}=\inf\left\{\lambda>0\::\:\int_{\Omega}\psi\left(\frac{|f(t)|}{\lambda}\right)dt\leq1\right\}<+\infty.
$$

\par Recall that a nonnegative function $\varphi$ defined on $[0;+\infty)$ is called quasiconcave if it satisfies the following conditions: $\varphi(0)=0$, $\varphi(t)$ is increasing, and $\varphi(t)/t$ is decreasing.

\par The Marcinkiewicz space $M_\varphi$ is defined as the set of all $f\in \mathcal{M}(\Omega)$ for which:
$$
||f||_{M_\varphi}=\sup_{0<t}\frac{1}{\varphi(t)}\int_0^tf^*(u)du<+\infty.
$$
\par It is worth noting that $(M_\varphi)_A = (M_\varphi)_B$, and $(M_\varphi)_A$ can be characterized as the set of functions $f\in\mathcal{M}$ (see \cite{kps}) that satisfy:
\begin{equation}
\label{estim_limit_decreasing_rearrangement_over_phi}
\lim_{t\to0+}\frac{1}{\varphi(t)}\int_0^tf^*(u)du=0.
\end{equation}

\par Additionally, when $\psi(t)=e^t-1$ and $\varphi(t)=t\ln(e/t)$, the corresponding Orlicz and Marcinkiewicz spaces coincide (see \cite{BS}), and we denote them as $e^L$ and $M_{\ln}$. Furthermore, it can be observed that (see \cite[Corollary 3.4.28]{EE}):
\begin{equation}
\label{estim_norm_equivalences}
||f||_{e^L}\asymp||f||_{M_{\ln}}\asymp\sup_{0<t\leq1}\frac{f^*(t)}{\ln(e/t)}.
\end{equation}

The following result was initially established in \cite{EGK} for the single-variable case, and our goal now is to extend it to the multi-variable scenario. Since the proof of this statement can be easily derived from the one provided in \cite{EGK}, we will omit it here.
\begin{thm}
\label{thm_necess_and_suff_C_closedc_in_X}
Let $X$ be a BFS on $\Omega$. The space $C(\Omega)$ of continuous functions is a closed linear subspace of $X$ if and only if there exists a positive constant $c$ such that for every rectangle $I \subset \Omega$, we have
\label{estim_C_is_closed_in_X}
$$c\leq ||\chi_{I}||_X.$$
\end{thm}

\section{Proof of Theorem \ref{thm_necess_suff_on_decreasing_rearrangement}}
The forthcoming proof closely follows the framework presented in \cite{kop-zv}. However, since we encounter some differences when extending the proof from the one-dimensional case to multiple dimensions, we have chosen to provide the complete proof for the sake of clarity.
\par Necessity. Since the space $C(\Omega)$ is closed in $L^{p(\cdot)}(\Omega)$, then by Theorem \ref{thm_necess_and_suff_C_closedc_in_X} there exists positive constant $d$ such that $d\leq ||\chi_{I}||_{p(\cdot)}$ for all rectangles $I$.
this implies $X_A\neq X_B$. Then by Proposition \ref{prop_X_A_equals_X_B} there exists $c>1$ such that 
\begin{equation}
\label{estim_int_c_pow_p_rearrangement_infty}
\int_0^1c^{p^*(t)}dt=+\infty.
\end{equation}
\par Consider two cases: 
\par Case 1) $p^*(\cdot)\in e^L$. Since (\ref{estim_int_c_pow_p_rearrangement_infty}) holds then function $p^*(\cdot)$ does not have absolute continuous norm that is $p^*(\cdot)\in e^L\backslash\left(e^L\right)_A$. Then by (\ref{estim_norm_equivalences}) we get that $p^*(\cdot)\in M_{\ln}\backslash(M_{\ln})_A$ and by (\ref{estim_limit_decreasing_rearrangement_over_phi}) it is obvious that
$$
\limsup_{t\to0+}\frac{1}{t\ln(e/t)}\cdot\int_0^tp^*(u)du>0,
$$
finally using ones more (\ref{estim_norm_equivalences}) from the last estimation we get (\ref{cond_necess_suffic}).

\par Case 2) $p^*(\cdot)\notin e^L$. Then by (\ref{estim_norm_equivalences})
$$
\sup_{0<t\leq1}\frac{p^*(t)}{\ln(e/t)}=+\infty,
$$
consequently (\ref{cond_necess_suffic}) holds.
The necessity part of the theorem proved.
\par Sufficiency. Let (\ref{cond_necess_suffic}) holds. For all $t\in[0;1]$ define function $h(t)=\min\{p^*(t), \ln(e/t)\}$. It is obvious that in this case holds
$$
\limsup_{t\to0+}\frac{h(t)}{\ln(e/t)}>0,
$$
then there exists a sequence $t_k\downarrow0$, such that
$$\frac{h(t_k)}{\ln(e/t_k)}\geq d,\quad k\in\mathbb{N},$$
for some positive number $d$. Now choose subsequence $(t_{k_n})$ such that $2t_{k_{n+1}}<t_{k_n}$, for all natural $n$. Since $t_k\downarrow0$, we can always choose such subsequence, so without loss of generality we can assume that sequence $(t_k)$ is already such.
\par Let given function $f$ defined by
$$
f(t)=d\cdot\ln(e/t_k), \:\:\: t\in(t_{k+1};t_k], \:\:\: k\in \mathbb{N} \quad\textnormal{and} \quad f(t)=1, \:\:\:t\in(t_1;1].
$$
It is clear that $h(t)\geq f(t)$ for all $t\in[0;1]$. Now choose positive number $c$ such that $c>e^{1/d}$ then we get 
\begin{equation}
\label{estim_c_pow_p_rearrangement_infty}
\int_0^1 c^{h(t)}dt=+\infty.
\end{equation}
Indeed,
$$
\int_0^1c^{h(t)}dt\geq\int_0^1c^{f(t)}dt>\int_{t_{k+1}}^{t_k}c^{d\cdot\ln(e/t_k)}dt=
$$
$$
=(t_k-t_{k+1})\cdot e^{d\cdot\ln c\cdot\ln(e/t_k)}>\frac{t_k}{2}\cdot\left(\frac{e}{t_k}\right)^{d\cdot\ln c}\to+\infty,\quad k\to+\infty.
$$
Choose decreasing sequence $\{a_k\}_{k\in\mathbb{N}}$, such that 
$$
\int_{a_{k+1}}^{a_k}c^{h(t)}dt=1.
$$
By (\ref{estim_c_pow_p_rearrangement_infty}) such sequence always can be chosen. Now let $\Delta_k=[a_{k+1};a_k]$, and $\{r_k\::\:k\in\mathbb{N}\}$ is a countable dense set in $[0;1]$. Define $b_k=-a_{k+1}+r_k$. Now let $A_k:=\Delta_k+b_k=[r_k;r_k+a_k-a_{k+1}]$. Let $g_k(t)=h(t) \cdot \chi_{\Delta_k}(t)$, $k\in\mathbb{N}$. Define functions $p_k(t)$ by the induction:
$$
p_1(t)=g_1(t-b_1)\chi_{[0;1]}(t),
$$
$$
p_k(t)=\left(p_{k-1}(t)(1-\chi_{\Delta_k}(t-b_k))+g_k(t-b_k)\right)\cdot\chi_{[0;1]}(t),\quad k>1.
$$

It is clear that $h(t)$ is decreasing and therefore $p_k(t)\leq p_{k+1}(t)$, for all $t\in[0;1]$ and all $k\in\mathbb{N}$. Also for all $k\in\mathbb{N}$ we have
\begin{equation}
\label{estim_int_p_k}
\int_0^1p_k(t)dt\leq\int_0^1h(t)dt\leq\int_0^1\ln(e/t)dt = 2.
\end{equation}
Now define $q(\cdot)$ function by
$$
q(t)=\lim_{k\to+\infty}p_k(t),\quad t\in[0;1].
$$
By (\ref{estim_int_p_k}) we get that the function $q(\cdot)$ is almost everywhere finite.
By the construction it is clear that $q^*(t)\leq h(t)\leq p^*(t)$. Now by the well known result (see \cite[Theorem 7.5]{BS}) there exists measure preserving transformation $\omega:[0;1]\to[0;1]$ such that $q(t)=q^*(\omega(t))$. Now define $\hat{p}(\cdot)$ by $\hat{p}(t) := p^{*}(\omega(t))$, $t\in[0;1]$. Since $q^*(t)\leq p^*(t)$ it is obvious that $q^*(\omega(t))\leq p^*(\omega(t))$, then for all $t\in(0;1)$ we get following inequality
\begin{equation}
\label{estim_bar_p_by_tilde_p}
q(t)\leq \hat{p}(t).
\end{equation}
Now, construct an exponential function $\bar{p}:\Omega \to [1,\infty)$, for which  the space of continuous functions will be a closed subspace inside its corresponding variable exponent Lebesgue space. For this, let's define measure-preserving mapping $\rho:\Omega \to [0;1]$, with the following rule: Suppose that, $x=(x_1,...,x_n)\in\Omega$ and for every $i\in\{1,...,n\}$ index, the representation of its corresponding coordinates be the following: $x_i=0.a_{i1}a_{i2}a_{i3}...$, then 
$$
\rho(x)=0.a_{11}a_{21}...a_{n1}a_{12}a_{22}...a_{n2}...
\:.$$ 
This mapping, which we have mentioned above, is well-known from literature. Thus, we can define the function $\bar{p}(x)=\hat{p}(\rho(x))$. In order to complete the proof, we should verify that the space of continuous function will be a closed subspace in its corresponding $L^{\bar{p}(\cdot)}(\Omega)$ space. For this purpose, we should show that there  exists a positive number $K$ such that for every rectangle $I \in \Omega $, we have $||\chi_{I}||_{\bar{p}}\ge K$. Consider any number $c > 1$, as in view of the fact that the set of binary rational numbers is dense everywhere in the set of all real numbers, for this reason we can find an $n$-dimensional binary rectangle $I^d$   for this given $n$-dimensional rectangle such that, $I^d \subset I $. Then because of the properties of function $\rho$, $c > 1$, and by \eqref{estim_bar_p_by_tilde_p}, we get 
$$
\int_Ic^{\bar{p}(x)}dx\ge \int_{I^d}c^{\bar{p}(x)}dx = \int_{(I^d)'}c^{\hat{p}(t)}dt\ge\int_{(I^d)'}c^{q(t)}dt,
$$ 
where $(I^d)'$ denotes one-dimensional binary interval taken out from $[0;1]$, for which $\rho(I^d)=(I^d)'$.  
By the construction of $q(\cdot)$ there exists number $k_0$ such that $E_{k_0}\subset (I^{d})'$. Then we get 
$$
\int_{(I^{d})'} c^{\hat{q}(t)} dt \ge \int_{E_{k_0}} c^{\hat{q}(t)} dt \ge \int_{E_{k_0}} c^{q_{k_0}(t)} dt = 
$$
$$ =\int_{E_{k_0}} c^{g_{k_0}(t-d_{k_0}} dt = \int_{r_{k_0}}^{r_{k_0} +t_{k_0} - t_{k_0 + 1}} c^{h(t-d_{k_0})\cdot\chi_{\Delta_{k_0}}(t-d_{k_0})} dt =
$$
$$
= \int_{t_{k_0 + 1}}^{t_{k_0}} c^{h(t)} dt \geq 1.
$$

Now by the definition of the norm in variable Lebesgue space and by the above estimations we get that for all $n$-dimensional rectangles $I \subset \Omega$ we have $||\chi_I||_{\bar{p}}>1/c$. By the Theorem \ref{thm_necess_and_suff_C_closedc_in_X} we get the proof of sufficiency of the Theorem \ref{thm_necess_suff_on_decreasing_rearrangement}.

% ----------------------------------------------------------------


\begin{thebibliography}{99}
\bibitem{BS} Bennet, C., Sharpley, R., Interpolation of operators, Pure Appl. Math. 129, Academic Press, 1988.

\bibitem{CUF}D. Cruz-Uribe, A. Fiorenza,
Variable Lebesgue Spaces: Foundations and Harmonic Analysis,
Birkh\"auser, Basel (2013).

\bibitem{DHHR} L. Diening, P. H\"ast\"o, P. Harjulehto and M. R$\stackrel{\circ}{\mbox{u}}$\v{z}i\v{c}ka, Lebesgue and
Sobolev spaces with variable exponents, Springer Lecture Notes, vol.
2017, Springer-Verlag, Berlin 2011.

\bibitem{EE} Edmunds, D., Evans, D., Hardy Operators, Function Spaces and Embeddings, Springer Berlin Heidelberg 2004.

\bibitem{EGK} Edmunds, D., Gogatishvili, A. and Kopaliani, T., Construction of function spaces close to $L^\infty$ with associate space close to $L^1$, J Fourier Anal Appl (2017). https://doi.org/10.1007/s00041-017-9574-2. 

\bibitem{ELN} Edmunds, D. E., Lang, J. and Nekvinda, A., On $L^{p(x)}$ norms, Proc. Roy. Soc. London, Series A 455 (1999), 219-225.



\bibitem{kop-zv} Kopaliani, T., Zviadadze Sh., Note on the variable exponent Lebesgue function spaces
close to $L^\infty$, J. Math. Anal. Appl., 474 (2019), 1463--1469.

\bibitem{kps} Krein, S. G., Petunin, Yu. I., Semenov, E. M., Interpolation of Linear Operators (American Mathematical Society, Providence, 1982). 


\bibitem{LN} Lang, J. and Nekvinda, A., A difference between continuous and absolutely
continuous norms in Banach function spaces, Czech. Math. J. 47 (2) (1997),
221-232.

\end{thebibliography}
\end{document}